\input amstex.tex
\documentstyle{amsppt}
%
%
%
%
%
%
%
%
%
%
\def\firstpage{1}\def\lastpage{1000}
\nopagenumbers

\hoffset=1.66cm\voffset=2.5cm 

\hsize=12.5cm\vsize=19.5cm

\font\eightrm=cmr8
\font\caps=cmcsc10                    
\font\Caps=cccsc10 scaled \magstep1   
\font\scaps=cmcsc8

%
\pageno=\firstpage
\def\folio{\rm\number\pageno}\output={\plainoutput}
\def\headfoot{25pt}
\def\makefootline{\baselineskip=\headfoot\line{\the\footline}}
\footline{\hfill\scaps Documenta Mathematica
     $\cdot$ Extra Volume ICM  1998  $\cdot$
    {\ifnum\pageno>\lastpage\else\number\firstpage--\lastpage\fi}\hfill}
\def\DocMath{{\def\th{\thinspace}\scaps Doc.\th Math.\th J.\th DMV}}
\def\makeheadline{
    \vbox to 0pt{\vskip-\headfoot\line{\vbox to8.5pt{}\the\headline}\vss}
    \nointerlineskip}
\headline={\ifnum\pageno=\firstpage{\DocMath\hfill\llap{\folio}}%
           \else{\ifodd\pageno\rightheadline\else\leftheadline\fi}\fi}
\def\rightheadline{\caps \hfill \rightheadtext    \hfill \llap{\folio}}
\def\leftheadline {\caps \rlap{\folio} \hfill \leftheadtext \hfill }
\def\leftheadtext{\ifnum\pageno>\lastpage\else\SAuthor\fi}
\def\rightheadtext{\STitle}
\def\TSkip{\bigskip}
\newbox\TheTitle{\obeylines\gdef\GetTitle #1
\ShortTitle  #2
\SubTitle    #3
\Author      #4
\ShortAuthor #5
\EndTitle
{\setbox\TheTitle=\vbox{\baselineskip=20pt\let\par=\cr\obeylines%
\halign{\centerline{\Caps##}\cr\noalign{\medskip}\cr#1\cr}}%
	\copy\TheTitle\TSkip\TSkip%
\def\next{#2}\ifx\next\empty\gdef\STitle{#1}\else\gdef\STitle{#2}\fi%
\def\next{#3}\ifx\next\empty%
    \else\setbox\TheTitle=\vbox{\baselineskip=20pt\let\par=\cr\obeylines%
    \halign{\centerline{\caps##} #3\cr}}\copy\TheTitle\TSkip\TSkip\fi%
\centerline{\caps #4}\TSkip\TSkip%
\def\next{#5}\ifx\next\empty\gdef\SAuthor{#4}\else\gdef\SAuthor{#5}\fi%
\catcode'015=5}}\def\Title{\obeylines\GetTitle}
\def\Abstract{\begingroup\narrower
    \parskip=\medskipamount\parindent=0pt{\caps Abstract. }}
\def\EndAbstract{\par\endgroup\TSkip}

\long\def\MSC#1\EndMSC{\def\arg{#1}\ifx\arg\empty\relax\else
     {\par\narrower\noindent%
     1991 Mathematics Subject Classification: #1\par}\fi}

\long\def\KEY#1\EndKEY{\def\arg{#1}\ifx\arg\empty\relax\else
	{\par\narrower\noindent Keywords and Phrases: #1\par}\fi\TSkip}

\newbox\TheAdd\def\Addresses{\vfill\copy\TheAdd\vfill
    \ifodd\number\lastpage\vfill\eject\phantom{.}\vfill\eject\fi}
{\obeylines\gdef\GetAddress #1
\Address #2
\Address #3
\Address #4
\EndAddress
{\def\xs{6truecm}
\setbox0=\vtop{{\obeylines\hsize=\xs#1\par}}\def\next{#2}
\ifx\next\empty 
     \setbox\TheAdd=\hbox to\hsize{\hfill\copy0\hfill}
\else\setbox1=\vtop{{\obeylines\hsize=\xs#2\par}}\def\next{#3}
\ifx\next\empty 
     \setbox\TheAdd=\hbox to\hsize{\hfill\copy0\hfill\copy1\hfill}
\else\setbox2=\vtop{{\obeylines\hsize=\xs#3\par}}\def\next{#4}
\ifx\next\empty\ 
     \setbox\TheAdd=\vtop{\hbox to\hsize{\hfill\copy0\hfill\copy1\hfill}
                \vskip20pt\hbox to\hsize{\hfill\copy2\hfill}}
\else\setbox3=\vtop{{\obeylines\hsize=\xs#4\par}}
     \setbox\TheAdd=\vtop{\hbox to\hsize{\hfill\copy0\hfill\copy1\hfill}
	        \vskip20pt\hbox to\hsize{\hfill\copy2\hfill\copy3\hfill}}
\fi\fi\fi\catcode'015=5}}\gdef\Address{\obeylines\GetAddress}

\hfuzz=0.1pt\tolerance=2000\emergencystretch=20pt\overfullrule=5pt
\Title
From Double Hecke Algebra to Analysis
\ShortTitle 
\SubTitle   
\Author 
Ivan Cherednik
\ShortAuthor 
\EndTitle
\Abstract 
We discuss $q$-counterparts of the Gauss integrals,
a new type of Gauss-Selberg sums at roots of 
unity, and  $q$-deformations 
of Riemann's zeta.
The paper contains general results, one-dimensional formulas,
and remarks about the current projects involving the double 
affine Hecke algebras.
\EndAbstract
\MSC 
\EndMSC
\KEY 
Hecke algebra, Fourier transform, spherical function,
Macdonald polynomial, Gauss integral, Gaussian sum,
metaplectic representation, Verlinde algebra, braid group, zeta function.
\EndKEY
\Address 
Ivan V. Cherednik
Dept. Math. UNC at Chapel Hill, 
Chapel Hill, NC 27599-3250, USA
chered\@math.unc.edu
\Address
\Address
\Address
\EndAddress

%
%
%
%
%

\font\eightrm=cmr8

\def\for{\  \hbox{ for } \ }
\def\if{ \ \hbox{ if } \ }

\def\where{\  \hbox{ where } \ }
\def\and{\  \hbox{ and } \ }

\def\equal{\buildrel  def \over =}

\def\om{\omega}

\def\th{\theta}
\def\al{\alpha}

\def\ga{\gamma}
\def\ep{\epsilon}

\def\ze{\zeta}


\def\vph{\varphi}

\def\tga{\tilde{\gamma}}

\def\tw{\tilde w}
\def\tW{\widetilde W}

\def\tz{\tilde z}
\def\tb{\tilde b}

\def\hW{\widehat{W}}

\def\hw{\hat{w}}

\def\hv{\hat{v}}

\def\C{\bold{C}}
\def\Q{\bold{Q}}

\def\R{\bold{R}}
\def\N{\bold{N}}
\def\Z{\bold{Z}}

\def\0{\bold{0}}


\font\germ=eufb10 
\def\goth#1{\hbox{\germ #1}}

\def\goth#1{\hbox{\germ #1}}
\def\equal{\buildrel  def \over =}

\def\HH{\hbox{${\Cal H}$\kern-5.2pt${\Cal H}$}}
\newif\ifacks
\font\smm=msbm10 at 12pt 
\def\symbol#1{\hbox{\smm #1}}
\def\lsmash{{\symbol n}}





\def\bf{\caps}
{\bf Introduction.}
\vskip 0.3cm

\noindent
The note is about the role of double affine Hecke algebras
in the unification of the 
classical zonal and  $p$-adic spherical functions and the
corresponding Fourier transforms. The new theory contains one
more parameter $q$ and, what is important,
dramatically improves the properties of the
Fourier transform. In contrast to the real
and $p$-adic theories, the  $q$-transform is selfdual
and has practically all other important properties of the classical
Fourier transform. Here I will mainly discuss the 
Fourier-invariance of the Gaussian. 

There are various applications. In combinatorics, they are
via the Macdonald polynomials. As $q\to 1,$
we  complete the Harish-Chandra theory of the 
spherical transform. The limit $q\to \infty$ covers
the $p$-adic Iwahori-Matsumoto-Macdonald theory. When $q$ is a root of
unity, we  generalize the Verlinde algebras, directly related
to quantum groups and Kac-Moody algebras,  and come to a new class of
Gauss-Selberg sums. 

However  the main applications could be of more analytic
nature. The representation of the double affine Hecke algebra 
generated by the Gaussian and its Fourier transform can be described in
full detail. So the next step is  to examine the spaces generated by 
Gaussian-type functions. The Fourier transforms
of the simplest examples  lead to
$q$-deformations of the classical zeta and $L$-functions.

\vskip 0.2cm
Of course there are other projects involving the double Hecke
algebras. I will mention at least some of them. The 
following  is far from being
complete.

1) {\it Macdonald's $q$-conjectures} [M1,M2]. 
Namely, the norm, duality, and  evaluation
conjectures [C1,C2]. My proof of the norm-formula
is similar to that from  [O1] in the differential case
(the duality and evaluation conjectures
collapse as $q\to 1$).
I would add to this list the Pieri rules [C2].
As to the nonsymmetric Macdonald polynomials, see [O2,M3,C3]. 
See also [M3,DS,Sa] about the $C^\vee C$  
(the Koornwinder polynomials), and [I,M4,C4] about the Aomoto
conjecture.

2) {\it $K$-theoretic interpretation}. I mean the papers [KL1,KK] and
more recent [GG,GKV]. Presumably it can lead to the 
Langlands-type description of irreducible representations of
double Hecke algebras, but the answer can be more complicated
than in [KL1] (see also recent Lusztig's papers on the representations
of affine Hecke algebras with unequal parameters). 
The Fourier transform is  misty in this approach. Let me add
here the strong Macdonald conjecture (Hanlon). 

3) {\it Induced and spherical representations}.
The classification of the spherical representations is much
simpler, as well as the irreducibility of
the induced ones. I used the technique of intertwiners in [C4]
following a similar  theory for the affine Hecke algebras.
The nonsymmetric polynomials form the simplest spherical
representation. There must be connections with [HO1].
The intertwiners also serve as creation operators for the
nonsymmetric Macdonald polynomials (the case of $GL$ is due
to [KS]).

4) {\it Radial parts via Dunkl operators}. The main references are
[D1,H,C5]. In the latter it was observed that the trigonometric 
differential Dunkl operators form the degenerate (graded) affine
Hecke algebra [L] ([Dr] for $GL_n$). The difference, elliptic, and 
difference-elliptic generalizations were introduced in [C6,C7,C8].
The nonsymmetric Macdonald polynomials are eigenfunctions
of the difference Dunkl operators. 
The connections with the KZ-equation play an important role here. 
I mean Matsuo's and my theorems from [Ma,C5,C6]. 
See also  [C9].

5) {\it Harmonic analysis}. In the rational-differential setup,
the definition of the generalized Bessel functions is from
[O3], the corresponding generalized Hankel transform was considered in [D2,J]
(see also [He]). In contrast to the spherical transform,
it is selfdual, as well as the difference generalization from [C2,C10].
The Mehta-Macdonald conjecture, directly related to the transform
of the Gaussian, was checked in [M1,O1] in the differential case 
and generalized in [C10]. See [HO2,O2,C11] about applications to the
Harish-Chandra theory.

6) {\it Roots of unity}. The construction from [C2] generalizes and,
at the same time, simplifies the Verlinde algebras. The latter
are formed by the so-called reduced representations
of quantum groups at roots of unity.
Another interpretation  is via the Kac-Moody algebras [KL2]
(due to Finkelberg for roots of unity). 
A valuable feature is the projective action of $PSL(2,\bold{Z})$
(cf. [K, Theorem 13.8]).
In [C3] the nonsymmetric polynomials are considered, which establishes
connections with the metaplectic (Weil) representations at roots of unity.

7) {\it Braids}. Concerning  $PSL(2,\bold{Z}),$ it acts projectively
on the double Hecke algebra itself.
The best known explanation (and proof) is based on the 
interpretation of this algebra as a quotient of the group algebra
of the fundamental group of the elliptic configuration space [C6].
The calculation is mainly due to [B] in the  $GL$-case.
For arbitrary root systems, it is similar to that from [Le], 
but our configuration space is different. 
Switching to the roots of unity, there may be applications
to the framed links including the Reshetikhin-Turaev invariants.

8) {\it Duality}. The previous discussion  was
 about arbitrary root systems. In the case of $GL$,
the theorem from [VV] establishes the duality between the double Hecke
algebras and the $q$-toroidal (double Kac-Moody)
algebras. It  generalizes the classical
Schur-Weyl duality, Jimbo's $q$-duality, and the affine analogues
from [Dr,C12]. When the  center charge is nontrivial it 
explains the results from [STU], which were recently extended by Uglov 
to irreducible representations of the Kac-Moody ${gl}_N$ of 
arbitrary positive integral levels.

\vskip 0.2cm

Let me also mention the relations of the symmetric Macdonald
polynomials (mainly of the $GL$-type) to:
a) the spherical functions on $q$-symmetric 
spaces (Noumi and others), 
b) the interpolation  polynomials (Macdonald, Lassalle,
Knop and  Sahi, Okounkov and Olshanski),
c) the quantum  $gl_N$ (Etingof, Kirillov Jr.),
d) the KZB-equation (---, ---,  Felder, Varchenko).
There are   connections with 
the affine Hecke algebra technique in 
the classical theory of $GL_N$ and $S_n.$ I mean, for instance,  [C12],  
papers of Nazarov and  Lascoux, Leclerc, Thibon,
and recent results towards the Kazhdan-Lusztig polynomials.

The coefficients of the  symmetric $GL$-polynomials
have interesting combinatorial properties
(Macdonald, Stanley,  Garsia, Haiman, ...). These polynomials appeared
in Kadell's work. Their norms are due to Macdonald,
the evaluation and duality conjectures were checked by Koornwinder, the
Macdonald operators were introduced independently by
Ruijsenaars together with  elliptic deformations.

Quite a few constructions can be extended 
to arbitrary finite groups
generated by complex reflections. For instance, 
the Dunkl operators and the
KZ-connection exist in this generality (Dunkl, Opdam, Malle). One can
try the affine and even the hyperbolic groups (Saito's root systems).

\vskip 0.2cm

{\bf 1. One-dimensional formulas}.
\vskip 0.2cm

\noindent
The starting point of many mathematical and physical theories 
is the formula:
$$
\eqalignno{
&2\int_{0}^\infty e^{-x^2}x^{2k}\hbox{d}x\ =\ \Gamma (k+1/2), \
\Re k>-1/2. 
&(1)
}
$$
Let us give some examples.
\vskip 0.2cm

{(a)} Its  generalization to the Bessel functions, namely,
the invariance of the Gaussian $e^{-x^2}$
with respect to the Hankel transform, 
is a cornerstone of the Plancherel formula.

{(b)} The following ``perturbation'' for the same $\Re k>-1/2$ 
$$
\eqalignno{
&\goth{Z}(k)\equal2\int_{0}^{\infty} 
(e^{x^2}+1)^{-1}x^{2k}\hbox{\, d}x = 
(1-2^{1/2-k})\Gamma (k+1/2)\zeta (k+1/2) &(2)
}
$$
is fundamental in the analytic number theory. 

{(c)} The multi-dimensional extension  due to Mehta  with  
$\prod_{1\le i<j\le n}(x_i-x_j)^{2k}$ instead of $ x^{2k}$
gave birth to the theory of matrix models and the Macdonald
theory with various applications in mathematics and physics.  

{(d)} Switching to the roots of unity, the Gauss formula 
$$
\eqalignno{
&\sum_{m=0}^{2N-1} e^{{\pi m^2 \over 2N}i}\ =\ 
(1+i)\sqrt{N},\ \ N\in {\bold N}&(3)
}
$$
can be considered as a certain counterpart of (1) at $k=0$.

{(e)} Replacing $x^{2k}$ by $\sinh(x)^{2k},$ we come to
the theory of spherical and hypergeometric functions and
to the spherical Fourier transform. The spherical transform
of  the Gaussian plays an  important role in the harmonic 
analysis on symmetric spaces.

To employ modern mathematics at full potential,
we do need to go from Bessel to  hypergeometric functions.
In contrast to the former, the latter  
can be studied, interpreted and generalized
by a variety of methods in the range from 
representation theory and algebraic geometry to integrable models 
and string theory.
However the straightforward passage  $x^{2k}\to \sinh(x)^{2k}$ creates
problems. The spherical transform is not selfdual anymore,
the formula (1) has no $\sinh$-counterpart, and
the Gaussian looses its Fourier-invariance.

\vskip 0.2cm
{\bf Difference setup.}
It was demonstrated recently that these important
features of the classical Fourier transform are restored for the kernel
$$
\eqalignno{
&\delta_k(x;q)\equal\prod_{j=0}^\infty {(1-q^{j+2x})(1-q^{j-2x})\over
(1-q^{j+k+2x})(1-q^{j+k-2x})},\ 0<q<1, \ k\in \bold{C}.
&(4)
}
$$ 
Actually the selfduality of the corresponding transform
can be expected a priori because the Macdonald truncated theta-function 
$\delta$ is 
a unification of
$\sinh(x)^{2k}$ and the Harish-Chandra function ($A_1$) serving the
inverse spherical transform.

As to (1), setting $q=\exp(-1/a),\ a>0,$
$$
\eqalignno{
&(-i)\int_{-\infty i}^{\infty i}q^{-x^2}\delta_k\hbox{\, d}x=
{2\sqrt{a\pi}}\prod_{j=0}^\infty
{1-q^{j+k}\over
 1-q^{j+2k}},\ \ \Re k>0.
&(5)
}
$$
Here both sides are well-defined for all $k$ except for the poles 
but coincide only when $ \Re k>0$, worse than in (1). This can be fixed
as follows:  
$$
\eqalignno{
&(-i)\int_{1/4-\infty i}^{1/4+\infty i}q^{-x^2}\mu_k\hbox{\, d}x=
{\sqrt{a\pi}}\prod_{j=1}^\infty
{1-q^{j+k}\over
 1-q^{j+2k}},\ \Re k>-1/2 \hbox{\ \ for}&(6)\cr
&\mu_k(x;q)\equal \prod_{j=0}^\infty {(1-q^{j+2x})(1-q^{j+1-2x})\over
(1-q^{j+k+2x})(1-q^{j+k+1-2x})},\ 0<q<1, \ k\in \bold{C}.
&(7)
}
$$
The limit of (6)  multiplied by $(a/4)^{k-1/2}$
as $a\to \infty$ is (1) in the imaginary variant.

Once we managed $\Gamma ,$  it would be unexcusable
not to try  (cf. (2))
$$
\eqalignno{
&\goth{Z}_q(k)\equal (-i)\int_{1/4-\infty i}^{1/4+\infty i} 
(q^{x^2}+1)^{-1}\mu_k\hbox{\, d}x
\hbox{\ \ for\ } \Re k>-1/2. &(8)
}
$$
It has a meromorphic continuation to all $k$ periodic in the imaginary
direction.
The limit of $(a/4)^{k-1/2}\goth{Z}_q$ as $a\to \infty$
is $\goth{Z}$ for all $k$ except for the poles.
The analytic continuation is based on the shift operator technique.
It seems  that all zeros of  $\goth{Z}_q(k)$ for
$a>1, \Re k>-1/2$ are $q$-deformations of the 
zeros of $\goth{Z}(k).$ 

\vskip 0.2cm
{\bf Jackson and Gauss sums.}
A most promising feature of special  $q$-functions is 
a posiblility to replace the integrals by sums,
the Jackson integrals. 

Let $\int_\sharp$ be the integration for the 
 path which  begins at $z=\epsilon i+\infty$, moves
to the left till $\epsilon i$, then down
through  the origin to  $-\epsilon i$, and then
returns down the positive real axis to $-\epsilon i+\infty$
(for small $\epsilon$). Then for $|\Im k|<2\epsilon, \Re k>0,$
$$
\eqalignno{
&{1\over 2i}\int_{\sharp} q^{x^2} \delta_k\hbox{\, d}x\ =
-{a\pi\over 2}\prod_{j=0}^\infty
{(1-q^{j+k})(1-q^{j-k})\over
 (1-q^{j+2k})(1-q^{j+1})}\times\goth{g}_q^{\sharp},\cr
&\goth{g}_q^{\sharp}(k)\equal\sum_{j=0}^\infty q^{{(k-j)^2\over 4}}
{1-q^{j+k}\over
1-q^{k}} \prod_{l=1}^j
{1-q^{l+2k-1}\over
 1-q^{l}}\ =
&(9)\cr
& q^{{k^2\over 4}}\prod_{j=1}^\infty
{(1-q^{j/2})(1-q^{j+k})(1+q^{j/2-1/4+k/2})(1+q^{j/2-1/4-k/2})\over
(1-q^j)}.
}
$$
The sum $\goth{g}_q^{\sharp}$ 
is the Jackson integral for a special choice ($k/2$) of the starting
point. The convergence of the sum (9) is for all $k.$ 
Similarly,
$$
\eqalignno{
&\goth{Z}_q^{\sharp}(k)\equal
-{a\pi\over 2}\prod_{j=0}^\infty
{(1-q^{j+k})(1-q^{j-k})\over
 (1-q^{j+2k})(1-q^{j+1})}\times\goth{z}_q^{\sharp},\cr
&\goth{z}_q^{\sharp}(k)\ =\ \sum_{j=0}^\infty  
q^{-kj} (q^{-{(k+j)^2\over 4}}+1)^{-1}
{1-q^{j+k}\over
1-q^{k}} \prod_{l=1}^j
{1-q^{l+2k-1}\over
 1-q^{l}}.
&(10)
}
$$
For all $k$ apart from the poles,
$\lim_{a\to \infty}({a\over 4})^{k-1/2}\goth{Z}_{q}^{\sharp}(k)=
\sin(\pi k)\goth{Z}(k).$

Numerically, it is likely that
all zeros of $\goth{Z}_q^{\sharp}$ in the 
strip 
$$
\{0\le\Im k<\sqrt{2\pi a}-\epsilon,\ \Re k>-1/2\} \hbox{\ \ for\ } 
a>2/\pi  \hbox{\  and \ small\ } \epsilon
$$ 
are deformations of the 
classical ones.
Moreover there is a strong tendency
for the  deformations of the zeros of the $\zeta(k+1/2)$-factor
to go to the right (big $a$).
They are  not expected in
the left half-plane before  $k= 1977.2714i$ (see [C13]).
 
When $q=\exp(2\pi i/N)$ and $k$ is a positive integer 
$\le N/2$ we come to the  Gauss-Selberg-type
sums:
$$
\eqalignno{
&\sum_{j=0}^{N-2k} q^{{(k-j)^2\over 4}}
{1-q^{j+k}\over
1-q^{k}} \prod_{l=1}^j
{1-q^{l+2k-1}\over
 1-q^{l}}=
\prod_{j=1}^k
(1-q^{j})^{-1}\sum_{m=0}^{2N-1} q^{m^2/4}.&(11)
}
$$
They resemble, for instance, [E,(1.2b)].
Substituting $k=[N/2]$ we  arrive at (3).

\vskip 0.2cm

\def\HH{\hbox{${\Cal H}$\kern-5.2pt${\Cal H}$}}

{\bf Double Hecke algebras} provide  justifications and
generalizations. In the $A_1$-case, 
$\HH\equal \C[\Cal{B}_q]/((T-t^{1/2})(T+t^{-1/2}))$ for
the group algebra of the group $\Cal{B}_q$ generated by
$T,X,Y,q^{1/4}$ with the relations
$$
\eqalignno{
&TXT=X^{-1},\ T^{-1}YT^{-1}=Y^{-1},\ Y^{-1}X^{-1}YXT^2=q^{-1/2}
&(12)
}
$$
for central $q^{1/4},t^{1/2}$. 
Renormalizing $T\to q^{-1/4}T,\ X\to q^{1/4}X,\
 Y\to q^{-1/4}Y,$
$$
\eqalignno{
&\Cal{B}_q\cong \Cal{B}_1\mod q^{1/4},\  \Cal{B}_1\cong
\pi_1(\{ E\times E\setminus \hbox{diag}\}/\bold{S}_2),
\ E=\hbox{elliptic \ curve,}
&(13)
}
$$
a special case of the calculation from [B].
The $T$ is the half-turn about the diagonal,
$X,Y$ correspond to  the  ``periods'' of $E.$

Thanks to the topological interpretation, the central
extension  $PSL^c_2(\bold{ Z})$ of
$PSL_2(\bold{ Z})$ (Steinberg)
acts on $\Cal{B}_1$ and $\HH$.
The automorphisms corresponding to
the generators ${11\choose 01},\
{10\choose 11}$
are as follows: 
$$
\eqalignno{
&\tau_+: Y\to q^{-1/4}XY,\ X\to X,\ \ \tau_-:  X\to q^{1/4}YX,\ Y\to Y,
&(14)
}
$$
fixing $T,q,t.$
When $t=1$ we get the well-known action of  
$SL_2(\bold{ Z})$ on the Weyl and  Heisenberg algebras (the latter as
$q\to 1$). Formally,  $\tau_+$ is the conjugation by
$q^{x^2}$ for $X$ represented here and later in the form $X= q^x.$  

The Mac\-do\-nald non\-sym\-met\-ric 
polynomials are eigenfunctions of  $Y$ in the following 
$\HH$-representation
in the space $\Cal{P}$ of the  Laurent polynomials of $q^x:$
$$
\eqalignno{
&T\to t^{1/2}s+(q^{2x}-1)^{-1}(t^{1/2}-t^{-1/2})(s-1), \ Y\to spT
&(15)
}
$$
for the reflection $sf(x)=f(-x)$ and the translation
$pf(x)=f(x+1/2).$ It is nothing else but the representation
of $\HH$ induced from the character $\chi(T)=t^{1/2}=\chi(Y).$
The Fourier transform (on the generalized functions) is associated
with the anti-involution $\{\vph: X\to Y^{-1}\to X\}$ of
$\HH$ preserving
$T,t,q.$ 

Combining $\tau_+$ and $\vph,$ we prove that
the Macdonald polynomials multiplied by  $q^{-x^2}$ are
eigenfunctions of the $q$-Fourier transform and get (6) for
$t=q^k.$

\vskip 0.2cm
When $q,k$ are from (11), let $q^x(m/2)=q^{m/2}$ for 
$m\in \Z,\ -N< m\le N,$
and
$$
\bowtie\ \equal \{m\mid \mu_k(m/2)\neq 0\} =
 \{-N+k+1,\ldots, -k,\ k+1, \ldots,N-k\}.
$$
The space $V_k=\hbox{Funct}(\bowtie)$ has a unique structure of
an (irreducible) $\HH$-module making the evaluation map
$\Cal{P}\ni f\mapsto f(m/2)\in V_k$ a  $\HH$-homomorphism. 
Setting $V_k=V^+_{k}\oplus V^-_{k}$ where
$T=\pm t^{\pm 1/2}$ on $V^{\pm}_k,$
the dimensions for $k<N/2$ are  $2(N-2k)=$ $ (N-2k+1)+(N-2k-1)$.
The components $V^\pm_{k}$ are $PSL^c_2(\bold{ Z})$-invariant.
Calculating its  action in  $V^+_{k}$ (which is a subalgebra 
of $V_k$) we come to the formulas
from [Ki,C2,C3];  $\ V^-_{k}$ is 
$PSL^c_2(\bold{ Z})$-isomorphic to  $V^+_{k+1}$.
For $k=1$ it is the Verlinde algebra. Involving the
shift operator, we get (11). 

We note that $V_k$  may have applications to the arithmetic theory of
coverings of elliptic curves ramified at one point thanks to  (13).

\vskip 0.2cm
{\bf 2. General results.}
\vskip 0.2cm

\noindent
Let $R=\{\al\}   \subset \R^n$ be a root system of type $A,B,...,F,G$
with respect to a euclidean form $(z,z')$ on $\R^n \ni z,z'$,
$W$ the Weyl group  generated by  the reflections $s_\al,$
$\al_1,...,\al_n$  simple roots, 
$R_{+}$ the set of positive  roots, $\om_1,...,\om_n$
the fundamental weights,
$ Q=\oplus^n_{i=1}\Z \al_i\subset$ $ P=\oplus^n_{i=1}\Z \om_i.$
We will also use coroots $\al^\vee =2\al/(\al,\al)$
and the corresponding $Q^\vee.$ The form will be normalized
by the condition  $(\th,\th)=2$ for the maximal coroot  $\th \in R^\vee_+$.

The affine Weyl group $\tW$ acts on $\tz=[z,\ze]\in \R^n\times \R$
and is generated by $s_i=s_{\al_i}$ and 
$s_{0}(\tz)\ =\  \tz+(z,\th)\al_0,$ for  $\al_0=[-\th,1].$
Setting  $b(\tz)=[z,\ze-(z,b)]$ for $b\in P$,
$ \tW=W\lsmash Q\subset$ $ \hW\equal W\lsmash P.$ 
We call the latter  the {\it extended affine
Weyl group}. It is generated over $\tW$
by the group $\pi\in \Pi\cong P/Q$ such
that $\pi$ leave the set $\al_0,\al^\vee_1,\ldots,\al^\vee_n$
invariant.  

The length $l(\hw)$  of $\hw = \pi\tw \in W^b,\ \pi\in \Pi, \tw\in W^a$
is by definition the length of the reduced decomposition of
$ \tw$ in terms of the simple reflections 
$s_i, 0\le i\le n.$  Given $ b\in P$, there is a unique
decomposition 
$$
\eqalignno{
&b= \pi_b w_b \hbox{\ \ such\ that \ }  w_b\in W,
\ l(b)=l(\pi_b)+l(w_b) \hbox{\ and \ } l(w_b)=\max.
&(16)
}
$$
Then $\Pi=\{\pi_{\om_r}\}$ for the
minuscule  $\om_r$:  $(\om_r,\al^\vee)\le 1$ for all $\al\in R_+$.

\vskip 0.2cm
{\bf Double Hecke algebras.}
Let  $q_\al= q^{(\al,\al)/2}, t_\al=q_\al^{k_\al}$ for 
$\{k_\al\}$ such that  $k_{w(\al)}= k_\al$ (all $w$),\
$t_i=t_{\al_i}$, $t_0=t_{\th},$\
$\rho_k= (1/2)\sum_{\al\in R_+} k_\al \al,$
$$
\eqalignno{
& X_{\tb}\ =\ \prod_{i=1}^n X_i^{l_i} q^{l} 
\if \tb=[b,l],
\  b=\sum_{i=1}^n l_i \om_i\in P,\ l\in (P,P)=(1/p)\Z
}
$$
for $p\in \N.$
By $\C^{\pm}_{q,t}[X]$  we mean the algebra of
polynomials in terms of $X_i^{\pm 1}$ 
over the field  $\C_{q,t}$ of rational functions of
$q^{1/(2p)}, t_\al^{1/2}.$
We will also use the evaluation
$X_{b}(q^z)\equal q^{(b,z)}.$

The  {\it double  affine Hecke algebra} $\HH\ $
is generated over the field $ \C_{ q,t}$ by 
the elements $\{ T_j,\ 0\le j\le n\}$, 
pairwise commutative $\{X_i\},$
and the group $\Pi$ where the following relations are imposed:

(o)\ \  $ (T_j-t_j^{1/2})(T_j+t_j^{-1/2})\ =\ 0,\ \ 0\ \le\ j\ \le\ n$;

(i)\ \ \ $ T_iT_jT_i...\ =\ T_jT_iT_j...,\ \ m_{ij}$ factors on each side;

(ii)\ \   $ \pi T_i\pi^{-1} =T_j,\ \pi X_b \pi^{-1}=
X_{\pi(b)} \if \pi\in \Pi,\ \pi(\al^\vee_i)=\al^\vee_j$; 

(iii)\  $T_iX_b T_i\ =\ X_b X_{a_i}^{-1} \if (b,\al^\vee_i)=1,\
1 \le i\le  n$;

(iv)\  $T_0X_b T_0\ =\ X_{s_0(b)}\ =\ X_b X_{\th} q^{-1}
\if (b,\th)=-1$;

(v)\ \ $T_iX_b\ =\ X_b T_i$ if $(b,\al^\vee_i)=0 \for 0 \le i\le  n.$

Here $m_{ij}$ are from the corresponding Coxeter relations.
Given $\tw \in \tW, \ \pi\in \Pi,$ the product
$T_{\pi\tw}\equal \pi T_{i_1}\cdots T_{i_l}$, where 
$\tw=s_{i_1}\cdots s_{i_l}$, $l=l(\tw)$,
does not depend on the choice of the reduced decomposition.
In particular, we arrive at the pairwise 
commutative elements  
$$
\eqalignno{
& Y_{b}\ =\  \prod_{i=1}^n Y_i^{l_i} \if  
b=\sum_{i=1}^n l_i\om_i\in P,\where  
 Y_i\equal T_{\om_i},
&(17)
}
$$
satisfying the relations
$\ T^{-1}_iY_b T^{-1}_i=Y_b Y_{a_i}^{-1}\ $ if $\ (b,\al^\vee_i)=1,\ \ $
$T_iY_b=Y_b T_i\ $if  $\ (b,\al^\vee_i)=0$, $\  1 \le i\le  n.$

The Fourier transform is related to  the anti-involution of $\HH$
$$
\eqalignno{
& \vph:\  X_i \to Y_i^{-1},\ \  Y_i \to X_i^{-1},\  \ T_i \to T_i,\
 t\to t,\  q\to  q,\ 1\le i\le n.
&(18)}
$$
The ``unitary'' representations are defined for the anti-involution
$$
\eqalignno{
  & X_i^*\ =\  X_i^{-1},\   Y_i^*\ =\  Y_i^{-1},\  
 T_i^* \ =\  T_i^{-1}, \
 t\to t^{-1},\
 q\to  q^{-1},\ 0\le i\le n.
}
$$

The next two automorphisms induce  a projective
action of $PSL_2(\Z):$ 
$$
\eqalignno{
& \tau_+: \ X_b \to X_b,\ \ Y_r \to X_rY_r q^{-(\om_r,\om_r)/2},\ 
Y_\th \to X_0^{-1}T_0^{-2}Y_\th,      \cr
& \tau_-: \ Y_b \to Y_b,\ \ X_r \to Y_r X_r q^{(\om_r,\om_r)/2},\ 
X_\th \to T_0X_0Y_\th^{-1} T_0, &(19)
}
$$
where $b\in P,\ \om_r$ are minuscule, $X_0\ =   q X_\th^{-1}$.
Obviously $\tau_-\  =\  \vph \tau_+\vph.$ The projectivity means
that $\tau_+^{-1}\tau_-\tau_+^{-1} =$ $
\tau_-\tau_+^{-1}\tau_-$.

\vskip 0.2cm
{\bf Polynomial representation.} 
Let $\hw(X_{\tb})\ =\ X_{\hw(\tb)}$ for $\hw \in \hW.$
Combining  the action of the group $\Pi,$ the multiplication
by  $X_b,$ and 
the {\it De\-ma\-zure-Lus\-ztig operators}  
$$
\eqalignno{
&T_j\  = \  t_j ^{1/2} s_j\ +\ 
(t_j^{1/2}-t_j^{-1/2})(X_{\al_j}-1)^{-1}(s_j-1),
\ 0\le j\le n,
&(20)
}
$$
we get a  representation of $\HH$ 
in  $\C^{\pm}_{q,t}[X].$

The coefficient of $X^0=1$ ({\it the constant term})
of a polynomilal
$f\in  \C^{\pm}_{q,t}[X]$
will be denoted by $\langle  f \rangle$. Let
$$
\eqalignno{
&\mu\ =\ \prod_{a \in R_+^\vee}
\prod_{i=0}^\infty {(1-X_\al q_\al^{i}) (1-X_\al^{-1}q_\al^{i+1})
\over
(1-X_\al t_\al q_\al^{i}) (1-X_\al^{-1}t_\al q_\al^{i+1})}.
&(21)
}
$$
We will consider $\mu$ as a Laurent series with the coefficients in
$\C[t][[q]]$.
The form $\langle \mu_0 f {g}^*\rangle$ makes the polynomial
representation unitary for 
$$
 X_b^*\ =\  X_{-b},\ t^*\ =\ t^{-1},\ q^*\ =\ q^{-1},
\ \mu_0=\mu_0/\langle \mu \rangle=\mu_0^*.
$$

The {\it Macdonald nonsymmetric polynomials} 
$\{e_b,b\in P\}$ are eigenvectors of
the operators $\{L_f\equal f(Y_1,\cdots, Y_n), f\in \C^{\pm}_{q,t}[X]\}$:
$$
\eqalignno{
&L_{f}(e_b)\ =\ f(q^{-b_\sharp})e_b, \where
b_\sharp\equal b-w_b^{-1}(\rho_k) \for w_b \hbox{\ from\ (16).}
&(22)
}
$$
They  are pairwise orthogonal with respect to the above pairing
and form a basis in $\C^{\pm}_{q,t}[X].$
The normalization  $\ep_b\equal e_b/ e_b(q^{-\rho_k})$
is the most  convenient in the harmonic 
analysis. For instance, the duality relations become especially simple:
$\ep_b(q^{c_\sharp})=\ep_c(q^{b_\sharp})$ for all $b,c\in P.$ 
The next  formula establishes that   
 $\ep_c$ multiplied by the Gaussian are eigenfunctions of
the difference Fourier transform:
$$
\eqalignno{
&\langle
\ep_b \ep_c^* \tga^{-1}\mu
\rangle = 
q^{(b_\sharp,b_\sharp)/2+(c_\sharp,c_\sharp)/2 -(\rho_k,\rho_k)} 
\ep_c^*(q^{b_\sharp})\times\cr
&\prod_{\al\in R_+}\prod_{ j=1}^{\infty}{ 
1- q_\al^{(\rho_k,\al^\vee)+j}\over
1-t_\al q_\al^{(\rho_k,\al^\vee)+j} }\  \for 
\tga^{-1}\equal\sum_{b\in P} q^{(b,b)/2}X_b.
&(23)
}
$$
When $b=c=0$ we get (5). Indeed, the series  for $\tga^{-1}$
is nothing else but the expansion of  
 $\ga^{-1}$ for $\ga=q^{x^2/2},$ where we set
$X_b=q^{x_b},\ x^2 = \Sigma_{i=1}^n x_{\om_i}x_ {\al_i^\vee}.$

\vskip 0.2cm
{\bf Jackson and Gauss sums.}
We fix generic $\xi\in \C^n$ and set $\langle f\rangle_\xi\equal
|W|^{-1}\sum_{w\in W, b\in B}f(q^{w(\xi)+b}).$
Here $f$ is a Laurent polynomial or any function well-defined
on $\{q^{w(\xi)+b}\}$.  We assume that $|q|<1$. 
For instance, 
$\langle \ga\rangle_\xi= \tga^{-1}(q^\xi)q^{(\xi,\xi)/2}.$ 
It is convenient to switch to 
$\mu^\circ(X,t)\equal$  $\mu^{-1}(X,t^{-1})$. 
Given $b,c\in P$,
$$
\eqalignno{
&\langle
\ep_b\, \ep_c^*\,\ga{\mu^\circ}
\rangle_\xi=
 q^{-(b_\sharp,b_\sharp)/2-(c_\sharp,c_\sharp)/2 +(\rho_k,\rho_k)}
\ep_c(q^{b_\sharp})\times \cr  
&|W|^{-1}\langle \ga\rangle_\xi 
 \prod_{\al\in R_+}\prod_{ j=0}^{\infty}{ 
1- t_\al^{-1}q_\al^{-(\rho_k,\al^\vee)+j}\over
1- q_\al^{-(\rho_k,\al^\vee)+j} }. 
&(24)
}
$$
For $\xi=-\rho_k$, (24) generalizes (9).
If $k\in \Z_+,$ then $\mu^\circ=q^{\hbox{\eightrm const}}\mu \in
\C^{\pm}_{q}[X]$, 
the  product in (24) is understood as the limit and becomes finite.

The proof of this formula and the previous one is based on the
analysis of the anti-involution (18) in the corresponding representations
of $\HH.$
Here it is the representation in   
$\Cal{F}=\hbox{Funct}(\hW,\C_{q,t}(q^{(\om_i,\xi)})).$ For $a,b\in P,$
$w\in W,
\hv\in \hW,$ we set
$$
\eqalignno{
&X_a(bw)=X_a(q^{b+w(\xi)}),\ X_a g(bw)=(X_a g)(bw),
\ \hv (g)(bw)=g(\hv^{-1}bw) 
}
$$
for $g\in\Cal{F}.$
It provides the  action of $X, \Pi.$ The $T$ act as follows:
$$
\eqalignno{
T_i(g)(\hw)\ &=\ 
{t_i^{1/2}q^{(\al_i,b+w(\xi))} - t_i^{-1/2}\over
q^{(\al_i,b+w(\xi))}- 1  }\
g(s_i \hw)\cr
&-{t_i^{1/2}-t_i^{-1/2}\over
q^{(\al_i,b+w(\xi))}- 1  }\
g(\hw) \for 0\le i\le n,
&(25)
}
$$
The formulas are closely connected with (20): the
natural evaluation map from
$\C^{\pm}_{q,t}[X]$ to $\Cal{F}$ is a $\HH$-homomorphism.
The unitarity is for  $\langle\mu_1 f g^* \rangle_\xi,$ where
the values of  $\mu_1=\mu/\mu(q^\xi)=\mu^\circ_1$ at $\hw$ are
$*$-invariant   ($\xi^*=\xi).$

Dropping the $X$-action, we get a deformation
of the regular representation
of the affine Hecke algebra generated by $T,\Pi$. Indeed, taking $\xi$ 
from the dominant affine Weyl chamber, (25)
tend to the $p$-adic formulas  from [Mat] when  $q\to \infty$
and $t$ are powers of $p.$
For $\xi=-\rho_k,$ the image
of the restriction map from $\Cal{F}$ to functions on the set 
$\{\pi_b, b\in P\},$ which is a $\HH$-homomorphism, generalizes
the spherical part of the regular representation.
The limit to the Harish-Chandra theory is $q \to 1$
where $k$ is fixed (the root multiplicity). See [He,C11].
\vskip 0.2cm

Now $q$ will be a primitive $N$-th root of unity,
$P_N= P/(P\cap NQ^\vee);$ the evaluations of Laurent polynomials
are functions on this set. Let $\langle f\rangle_N\equal$
$\sum_{b\in P_N}f(q^{b}).$ 
We assume that $k_\al\in \Z_+$ for all $\al\in R$ and
 $\mu(q^{-\rho_k})\neq 0.$ We also pick $q$ to ensure the existence 
of the Gaussian: $q^{(b,b)/2}=1$ for all $b\in P\cap NQ^\vee.$
It means that when $N$ is odd and the root system is either
$B$ or $C_{4l+2}$ one takes  $q=\exp(4\pi i m/N)$ for $(m,N)=1, 0< 2m<N.$
Otherwise it is arbitrary.

We claim that the formula
(24) holds for  $\langle \ \rangle_N$ instead of 
 $\langle \ \rangle_\xi$ provided the existence of the
nonsymmetric polynomials. It readily gives (11) for $b=0=c.$ 

Given $b'\in P_N$ such that  $\mu(q^{b'})\neq 0,$  
at least one $\ep_{b}$ exists with $b_\sharp$ equal to 
$b'$ in $P_N.$
Denoting the set of all such $b'$
by $P'_N,$ the space $\hbox{Funct}(P'_N,\Q(q^{1/(2p)}))$
is an algebra and a $\HH$-module isomorphic to the quotient
of the polynomial representation by the radical of
the pairing  $\langle \mu f {g}^*\rangle_N.$ The radical also 
coincides with the set of  polynomials $f$ such that
$(g(Y)(f))(q^{-\rho_k})=0$ for all Laurent polynomials $g.$ 
The evaluations of $\ep_{b}$ depend only on the
images of  $b_\sharp$ in  $P_N$ and form a basis of this module.
The evaluations of the symmetric polynomials
constitute the {\it generalized Verlinde algebra}.

\vskip 0.2cm
{\it Acknowlegments.} {\eightrm 
Partially supported by NSF grant
DMS-9622829  and the Guggenheim Fellowship. The paper was completed
at the University Paris 7, the author is grateful for the kind invitation.}

\comment
\item{[N]}
Nagell, Trygve
Introduction to number theory
Chelsea Publishing Company, N.Y. 1981
pg 177, The Gaussian sum (53. Theorem 99, Chapter V
\endcomment

\vskip 0.2cm
\noindent
{\bf References.}
\vskip 0.2cm

\parindent=30pt
\eightrm
\baselineskip=10pt 
\font\ebf=cmbx8

\item{[B]}
Birman, J.: 
On braid groups. Communs on Pure and Appl. Math.  
{\ebf 22}, 41--72 (1969).

\item{[C1]}
Cherednik, I.:
Double affine Hecke algebras and Macdonald's conjectures.
 Annals of Math.
{\ebf 141}, 
 191--216 (1995).

\item{[C2]}
---:
 Macdonald's evaluation conjectures and difference Fourier
transform.
 Invent. Math. {\ebf 122}, 119--145
(1995).

\item{[C3]}
---:  Nonsymmetric Macdonald polynomialas.
IMRN {\ebf 10}, 483--515 (1995).

\item{[C4]}
---: Intertwining operators of double affine Hecke algebras.
Selecta Math. New s. {\ebf 3}, 459--495  (1997).

\item{[C5]}
---: Integration of Quantum many-body problems by affine
Knizhnik-Za\-mo\-lod\-chi\-kov equations. Advances in Math.  
{\ebf 106}, 
65--95  (1995).

\item{[C6]}
---: Double affine Hecke algebras, Knizhnik--Zamolodchikov
equations, and Macdonald's operators. IMRN (Duke Math. J.) {\ebf 9},
171--180 (1992).

\item{[C7]}
---:
Elliptic quantum many-body problem and double affine
 Knizhnik - Zamolodchikov equation.
Communs Math. Phys. {\ebf 169}:2, 441--461 (1995).

\item{[C8]}
---: Difference-elliptic operators and root systems, 
IMRN {\ebf 1}, 43--59 (1995).

\item{[C9]}
---:
Lectures on Knizhnik-Zamolodchikov
equations and  Hecke
algebras.  MSJ Memoirs  {\ebf 1}  (1998).

\item{[C10]}
---:
Difference Macdonald-Mehta conjectures.
IMRN {\ebf 10}, 449--467 (1997).

\item{[C11]}
---:
Inverse Harish--Chandra transform and difference operators.
IMRN {\ebf 15}, 733--750  (1997).

\item{[C12]}
---: A new interpretation of Gelfand-Tzetlin bases.
 Duke Math. J. {\ebf 54}:2,  563--577 (1987).

\item{[C13]}
---:
On $q$-analogues of Riemann's zeta. Preprint (1998).

\item{[DS]}
Diejen, J.F. van, Stockman, J.V.:
Multivariable $q$-Racah polynomials.
Duke Math. J. {\ebf 91}, 89--136 (1998).

\item{[Dr]}
Drinfeld, V.G.: Degenerate affine Hecke algebras and Yangians.
Funct. Anal. and Appl. {\ebf 20}, 69--70 (1986).

\item{[D1]}
Dunkl, C.F.: Diffirential-difference operators associated to
reflection groups. Trans. AMS. {\ebf 311}, 167--183 (1989).

\item{[D2]}
---:
Hankel transforms associated to finite reflection groups.
Contemp. Math. {\ebf 138},  123--138 (1992).

\item{[E]}
Evans, R.J.:
The evaluation of Selberg character sums.
L'Enseignement Math.   {\ebf 37}, 235--248 (1991).

\item{[GG]}
Garland, H., Grojnowski, I.:
Affine Hecke algebras associated to Kac-Moody groups.
Preprint (1995).

\item{[GKV]} 
Ginzburg, V., Kapranov, M., Vasserot, E.:
Residue construction of Hecke algebras. Preprint (1995).

\item{[H]}
Heckman, G.J.: An elementary approach to the hypergeometric shift 
operator of Opdam. Invent. math. {\ebf 103} 341--350 (1991).

\item{[HO1]}
Heckman, G.J., Opdam, E.M.:
Harmonic analysis for affine Hecke algebras.
Preprint (1996).

\item{[HO2]}
---: Root systems and hypergeometric functions I.
Comp. Math. {\ebf 64}, 329--352 (1987).

\item {[He]}
{Helgason, S.}:
Groups and geometric analysis. 
Academic Press, New York (1984).

\item {[I]}
Ito, M.: On a theta product formula
for Jackson integrals associated with root
systems of rank two.
Preprint (1996).

\item {[J]}
{Jeu, M.F.E. de}:
The Dunkl transform. Invent. Math. {\ebf 113},  147--162 (1993).

\item {[K]}
{Kac, V.G.}:
Infinite dimensional Lie algebras.
Cambridge University Press, Cambridge (1990).

\item {[Ki]}
Kirillov, A.  Jr.:
Inner product on conformal blocks and Macdonald's
polynomials at roots of unity. Preprint (1995).

\item {[KL1]}
Kazhdan, D., Lusztig, G.:
Proof of the Deligne-Langlands conjecture for Hecke algebras.
Invent. Math. {\ebf 87}, 153--215 (1987).

\item {[KL2]}
---: Tensor structures arising from affine Lie algebras. III.
J. of AMS {\ebf  7}, 335--381 (1994).

\item {[KS]}
Knop, F., Sahi, S.:
A recursion  and a combinatorial formula for Jack
polynomials,
Preprint (1996), to appear in Invent. Math.

\item {[KK]}
Kostant, B., Kumar, S.:
T-Equivariant K-theory of generalized flag varieties.
J. Diff. Geometry {\ebf 32},  549--603 (1990).

\item{[Le]}
Lek, H. van der: Extended Artin groups. Proc. Symp. Pure Math.
{\ebf 40}:2, 117--122 (1981).

\item{[L]}
Lusztig, G.: Affine Hecke algebras and their graded version.
J. of the AMS {\ebf 2}:3, 599--685 (1989).

\item{[M1]}
Macdonald, I.G.: Some conjectures for root systems. SIAM J.
Math. An. {\ebf 13}, 988--1007 (1982).

\item{[M2]}
---: Orthogonal polynomials associated with root
systems, Preprint(1988).

\item{[M3]}
---: Affine Hecke algebras and orthogonal polynomials.
S\'eminaire Bourbaki {\ebf  47}:797, 01--18 (1995).

\item{[M4]}
---: A formal identity for affine root systems,
Preprint (1996).

\item{[Ma]}
Matsuo, A.: Knizhnik-Zamolodchikov type equations and zonal
spherical functions.
Invent. Math. {\ebf 110}, 95--121 (1992).

\item{[Mat]}
Matsumoto, H.: Analyse harmonique dans les systemes de
Tits bornologiques de type affine.
Lecture Notes in Math.  {\ebf 590} (1977).

\item{[O1]}
Opdam, E.M.:
Some applications of hypergeometric shift
operators. Invent. Math. {\ebf  98}, 1--18 (1989).

\item{[O2]}
---: Harmonic analysis for certain representations of
gra\-ded Hec\-ke al\-gebras. 
Acta Math. {\ebf 175}, 75--121 (1995).

\item{[O3]}
---:
Dunkl operators, Bessel functions and the discriminant of
a finite Coxeter group,
Comp. Math. {\ebf 85}, 333--373 (1993).

\item{[Sa]}
Sahi, S.: 
Nonsymmetric Koornwinder polynomials and duality.
Preprint (1996), to appear in Annals of Math.

\item{[STU]}
Saito, Y., Takemura,  K., Uglov, D.:
Toroidal actions on level-1 modules of $U_q(\hat{sl_n})$.
Transformation Groups {\ebf 3}, 75--102 (1998).

\item{[VV]} 
Varagnolo, M., Vasserot, E.:
Double-loop algebras and the Fock space.
Preprint (1996), to appear in  Invent.Math.

\Addresses
\end